\documentclass[11pt,fleqn,leqno]{article}
\usepackage{amsmath}
\usepackage{amsfonts}
\usepackage[dvips]{graphicx}
\usepackage{epsf}
\newtheorem{thm}{Theorem}[section]
\newcommand{\be}{\begin{equation}}
\newcommand{\ee}{\end{equation}}
\newcommand{\ba}{\begin{array}}
\newcommand{\ea}{\end{array}}

\renewcommand{\a}{\alpha}
\renewcommand{\b}{\beta}
\renewcommand{\l}{\lambda}
\renewcommand{\t}{\theta}

\newcommand{\bea}{\begin{eqnarray}}
\newcommand{\eea}{\end{eqnarray}}
\newcommand{\ga}{\gamma}
\newcommand{\w}{\omega}

\newcommand{\Ga}{\Gamma}
\newcommand{\la}{\lambda}

\newcommand{\ck}{{\cal K}}

\newcommand{\ei}{{\rm e}^{i\t}}
\newcommand{\emi}{{\rm e}^{-i\t}}
\renewcommand{\sf}{_{\infty}}
\newcommand{\dt}{{\frac{d\t}{2\pi}}}

\begin{document}
\title{Small eigenvalues of large Hankel matrices: The indeterminate case
\thanks{This research is partially supported by the EPSRC GR/M16580
 and NSF grant DMS 99-70865}}
\author{Christian Berg, Yang Chen and Mourad E. H. Ismail}
\maketitle
                                                   
\date{}

\begin{abstract} In this paper we characterise the indeterminate case
by the eigenvalues of the Hankel matrices being bounded below by a strictly
positive constant. An explicit lower bound is given in terms of the orthonormal
polynomials and we find expressions for this lower bound in a number of
indeterminate moment problems.
\end{abstract}

\bigskip

\setcounter{section}{0}
\setcounter{equation}{0}
\setcounter{thm}{0}

\section{Introduction}
Let $\a$ be a positive measure on {\bf R} with infinite support
and finite  moments of all orders
\bea
s_n=s_n(\a)=\int_{\bf R} x^n d\a(x).
\eea
 With $\a$ we associate
the infinite Hankel matrix ${\cal H}_\infty = \{H_{jk}\}$,
\bea
H_{jk}= s_{j+k}.
\eea
Let ${\cal H}_N$ be the $(N+1)\times(N+1)$ matrix whose entries are $H_{jk},
0 \le j, k \le N$. Since ${\cal H}_N$ is positive definite, then all its
eigenvalues are positive. The large $N$ asymptotics of
the smallest eigenvalue, denoted as $\la_N$, of the Hankel matrix
${\cal H}_N$ has been studied in papers by Szeg\"o \cite{Sz},
Widom and Wilf \cite{Wid}, Chen and
Lawrence \cite{Ch:La}. See also the monograph by Wilf \cite{Wi}.
All the cases considered
by these authors are determinate moment problems, and it was shown in
 each case that
$\l_N \to 0$, and asymptotic results were obtained about how fast $\la_N$
tends to zero.

The smallest eigenvalue can be obtained from the classical Rayleigh
quotient:
\bea
\la_N={\rm min}\left\{\sum_{j=0}^{N}\sum_{k=0}^{N}\;s_{j+k}v_jv_k:
\sum_{k=0}^{N}v_j^2=1, v_j\in {\bf R},\;0\leq j\leq N \right\}.
\eea
It follows  that $\la_N$ is a decreasing function of $N.$

The main result of this paper is Theorem 1.1, which we state next.

\begin{thm}
The moment problem associated with the moments (1.1) is determinate if
 and only if $\lim_{N\to\infty}\la_N=0$.
\end{thm}

We shall compare this result with a theorem of Hamburger \cite[Satz XXXI]{H},
cf. \cite[p.83]{Ak} or \cite[p.70]{Sh:Ta}

Let $\mu_N$ be the minimum of the Hankel form ${\cal H}_N$ on the hyper-plane
$v_0=1$, i.e.
\bea
\mu_N={\rm min}\left\{\sum_{j=0}^{N}\sum_{k=0}^{N}\;s_{j+k}v_jv_k:
v_0=1, v_j\in {\bf R},\;0\leq j\leq N \right\}.
\eea
and let $\mu_N'$ be the corresponding minimum for the moment sequence
$s_n'=s_{n+2}, n\geq 0$, i.e.
\bea
\mbox{}
\mu_N'&=&{\rm min}\left\{\sum_{j=0}^{N}\sum_{k=0}^{N}\;s_{j+k+2}v_j'v_k':
v_0'=1, v_j'\in {\bf R},\;0\leq j\leq N \right\}\nonumber\\
\mbox{}
&=&{\rm min}\left\{\sum_{j=0}^{N+1}\sum_{k=0}^{N+1}\;s_{j+k}v_jv_k:
v_0=0, v_1=1, v_j\in {\bf R},\;0\leq j\leq N+1 \right\}.\nonumber
\eea
The theorem of Hamburger can be stated that the moment problem is determinate
if and only if at least one of the limits $\lim_{N\to\infty}\mu_N$,
$\lim_{N\to\infty}\mu_N'$ are zero.

It is clear from (1.3), (1.4) that $\mu_N\geq\la_N$ and similarly
$\mu_N'\geq\la_{N+1}$. From these inequalities and Hamburger's theorem, we
obtain the \lq\lq only if\rq\rq\;\; statement in Theorem 1.1.
 The \lq\lq if\rq\rq\;\;
statement will be proved by finding a positive lower bound for
 the eigenvalues $\la_N$, cf. Theorem 1.2 below.

 We think that Theorem 1.1 has the advantage over the theorem of Hamburger
  that it involves only the moment sequence $(s_n)$ and not the shifted
  sequence $(s_{n+2})$. In section 2 we give another proof of the \lq\lq
  only if\rq\rq \;\;statement to make the proof of Theorem 1.1 independent of
   Hamburgers theorem.

If
\bea
\pi_{N}(x):=\sum_{j=0}^{N}v_jx^j,
\eea
then a simple calculation shows that
\bea
\sum_{0\leq j,\;k\leq N}s_{j+k}v_jv_k=\int_{E}\pi_N^2(x)\;d\a(x),
\eea
and
\bea
\sum_{k=0}^{N}v_k^2=\int_{0}^{2\pi}\left|\pi_{N}({\rm e}^{i\t})
\right|^2\;\frac{d\t}{2\pi}.
\eea
We could also study the reciprocal of $\la_N$ given by
\bea
\frac{1}{\la_N}={\rm max} \left\{\int_{0}^{2\pi}
\left|\pi_N({\rm e}^{i\t})\right|^2\;\frac{d\t}{2\pi}:\;\pi_N,\;
\int_{E}\pi_N^2(x) d\a(x)=1\right\}.
\eea
Let $\{p_k\}$ denote the orthonormal polynomials with respect to $\a$,
normalised so that $p_k$ has positive leading coefficient.

We recall that the moment problem is indeterminate, cf.
 \cite{Ak},\cite{Sh:Ta}, if and only
if there exists a non-real number $z_0$ such that
\bea
\sum_{k=0}^\infty |p_k(z_0)|^2<\infty.
\eea
In the indeterminate case
the series in (1.9) actually converges  for all $z_0$ in {\bf C}, uniformly on
 compact sets. In the determinate case the series in (1.9) diverges for all
 non-real $z_0$ and also for all real numbers except the at most countably many
 points, where $\a$ has a positive discrete mass.

 If we expand the polynomial (1.5)
as a linear combination of the  orthonormal system
\bea
\pi_{N}(x)=\sum_{j=0}^{N}c_jp_j(x),\nonumber
\eea
then
\bea
\int_{0}^{2\pi}\left|\pi_N ({\rm e}^{i\t})\right|^2\dt
=\sum_{0\leq j,\;k\leq N}c_jc_k\int_{0}^{2\pi}p_j ({\rm e}^{i\t})
p_k ({\rm e}^{-i\t})\dt
=\sum_{0\leq j,\;k\leq N}\ck_{jk}c_jc_k,\nonumber
\eea
where we have defined
\bea
\ck_{jk}=\int_{0}^{2\pi}p_j({\rm e}^{i\t})
p_k({\rm e}^{-i\t})\frac{d\t}{2\pi}.
\eea
Thus
\bea
\frac{1}{\la_N}={\rm max}
\left\{\sum_{0\leq j,k\leq N}\ck_{jk}c_jc_k:\;c_j,\;\sum_{j=0}^{N}c_j^2=1
\right\}.
\eea
Since the eigenvalues of the matrix $({\cal K}_{jk})_{0\le j,k\le N}$
 are positive, and their sum
is its trace, then
\bea
\frac{1}{\la_N}\leq \sum_{k=0}^{N}{\cal K}_{kk}=
\int_{0}^{2\pi}\sum_{k=0}^{N}\left|p_k({\rm e}^{i\t})
\right|^2\;\dt.
\eea

Thus in the case of indeterminacy,
\bea
\frac{1}{\la_N}\leq
\int_{0}^{2\pi}\sum_{k=0}^\infty\left|p_k({\rm e}^{i\t})
\right|^2\;\dt<\infty,
\eea
which shows that
\bea
\lim_{N\to\infty}\la_N\geq
\left(\int_{0}^{2\pi}\frac{1}{\rho(\ei)}\;\dt\right)^{-1},
\eea
where
\bea
\rho(\ei)
=\left(\sum_{k=0}^{\infty}\left|p_k({\rm e}^{i\t})\right|^2
\right)^{-1}.
\eea
We recall that for $z\in {\bf C}\setminus {\bf R}$ the number
 $\rho(z)/|z-{\overline z}|$ is the radius of the Weyl circle at $z$.

The above argument establishes the following result:

\begin{thm} In the indeterminate case the smallest eigenvalue $\la_N$ of
the Hankel matrix ${\cal H}_N$ is bounded below by the harmonic mean
 of the function $\rho$ along the unit circle.
\end{thm}

We shall conclude this paper with examples, where we have calculated
or estimated the quantity
\bea
\rho_0=\int_0^{2\pi}\sum_{k=0}^{\infty}\biggl|p_k\left({\rm e}^{i\t}\right)
\biggr|^2\;\dt.
\eea
This will be done for the moment problems associated with
 the Stieltjes-Wigert polynomials, cf. \cite{Ch},\cite{Sz2}, the
  Al-Salam-Carlitz polynomials \cite{Al:Ca}, the
symmetrized version of polynomials of Berg-Valent (\cite{Be:Va}) leading
 to a Freud-like weight \cite{Ch:Is}, and the
$q^{-1}$-Hermite  polynomials of Ismail and Masson \cite{Is:Ma}.

If we introduce the coefficients of the orthonormal polynomials as
\bea
p_k(x)=\sum_{j=0}^k\b_{k,j}x^j
\eea
then
\bea\int_0^{2\pi} |p_k(\ei)|^2\;\dt=\sum_{j=0}^k\b_{k,j}^2,
\nonumber
\eea
and therefore
\bea
\rho_0=\sum_{k=0}^\infty\sum_{j=0}^k \b_{k,j}^2.
\eea

Another possibility for calculating $\rho_0$ is to use the entire functions
 $B,D$ from the Nevanlinna matrix since
it is well known that \cite[p.123]{Ak}
\bea
\sum_{k=0}^\infty |p_k(z)|^2 =
\frac{B(z)D(\overline{z})- D(z)B(\overline{z})}{z - \overline{z}}.
\eea

It follows that
\bea
\sum_{k=0}^\infty |p_k(\ei)|^2 = {\rm Im}\{B(\ei)D(\emi)\}/\sin \t.
\eea

\bigskip

\setcounter{equation}{0}
\setcounter{thm}{0}

\section{Indeterminate Moment Problems}

In this section we shall give a proof of Theorem 1.1 which is independent
of Hamburgers result.
 We have already
established that if $\lim_{N\to\infty}\la_N=0$, then the problem is
determinate. We shall next prove that if $\la_N\geq \gamma$ for all $N$, where
$\gamma>0$, then the problem is indeterminate. Since $1/\la_N\leq 1/\gamma$
 for all $N$, and $1/\la_N$ is the biggest eigenvalue of the positive
 definite matrix $(\ck_{jk})_{0\leq j,k\leq N}$, we get
\bea
\sum_{0\leq j,k\leq N}\ck_{jk}c_j\overline{c_k}\le\frac{1}{\gamma}
\sum_{j=0}^{N}|c_j|^2,
\eea
for all vectors $(c_0,\ldots,c_N)\in{\bf C}^{N+1}$.
If we consider an arbitrary complex polynomial $p$ of degree $\le N$ written
as $p(x)=\sum_{k=0}^Nc_kp_k(x)$, the inequality (2.1) can be formulated
\bea
\int_0^{2\pi}|p(\ei)|^2\;\dt\leq\frac{1}{\gamma}
\int |p(x)|^2\;d\a(x).
\eea

Let now $z_0$ be an arbitrary non-real number in the open unit disc. By
the Cauchy integral formula
$$
p(z_0)=\frac{1}{2\pi}\int_0^{2\pi}\frac{p(\ei)}{\ei-z_0}\ei\;d\t,
$$
and therefore
\bea
|p(z_0)|^2\leq \int_0^{2\pi}|p(\ei)|^2\;\dt
\int_0^{2\pi}\frac{1}{|\ei-z_0|^2}\;\dt.
\eea

Combined with (2.2) we see that there is a constant $K$ such that for all
complex polynomials $p$
\bea
|p(z_0)|^2\leq K\int |p(x)|^2\;d\a(x),
\eea
where $K=1/(\gamma(1-|z_0|^2)).$

This inequality implies indeterminacy in the following way. Applying it to
 the polynomial
\bea
p(x)=\sum_{k=0}^N p_k(\overline{z_0})p_k(x),\nonumber
\eea
we get
\bea                                                 
\sum_{k=0}^N|p_k(z_0)|^2\leq K,
\eea
and since $N$ is arbitrary, indeterminacy follows.

{\bf Remark}. We see that the infinite positive definite matrix 
${\cal K}_\infty=\{{\cal K}_{j,k}\}$ is bounded on $\ell^2$ if and only if
$\la_N\geq\ga$ for all $N$ for some $\ga>0$. Furthermore
${\cal K}_\infty$ is of trace class if and only if $\rho_0<\infty$.
The result of Theorem 1.1 can be reformulated to say that boundedness implies
 trace class  for this family of operators.

\bigskip

\setcounter{equation}{0}
\setcounter{thm}{0}

\section{Examples}

We shall follow the notation and terminology  for $q$-special 
functions as those in Gasper and Rahman \cite{Ga:Ra}.

{\bf Example 1}. The Stieltjes-Wigert Polynomials.  

These polynomials are orthonormal with respect to the weight function
\bea
\omega(x)=\frac{k}{\sqrt{\pi}}\exp(-k^2(\log x)^2), \quad x > 0,
\eea
where $k>0$ is a positive parameter, cf. \cite{Ch},\cite{Sz2}. They
are given by
\bea
p_n(x)=(-1)^nq^{\frac{n}{2}+\frac14}(q;q)_n^{-\frac12}\sum_{k=0}^n
     \binom {n}{k}_q q^{k^2}(-q^{\frac12}x)^k,
\eea
where we have defined $q=\exp\{-(2k^2)^{-1}\}$.

It follows by (1.18) that
\bea
\mbox{}
\rho_0&=&\sum_{n=0}^\infty \frac{q^{n+\frac12}}{(q;q)_n}\sum_{k=0}^n
q^{k(2k+1)}\binom {n}{k}_q^2\\
\mbox{} &=& \sum_{k=0}^\infty q^{2k^2+k+\frac12}\sum_{n=k}^\infty
\frac{q^n}{(q;q)_n}\binom {n}{k}_q^2.
\nonumber
\eea
Putting $n=k+j,$ the inner sum is
\bea
\sum_{j=0}^\infty \frac{q^{k+j}}{(q;q)_k^2}\frac{(q;q)_{k+j}}{(q;q)_j^2}=
\frac{q^k}{(q;q)_k}{}_2\phi_1(q^{k+1},0;q;q,q)
\nonumber
\eea
and hence
\bea
\rho_0=\sum_{k=0}^\infty\frac{q^{2(k+\frac12)^2}}{(q;q)_k}
{}_2\phi_1(0, q^{k+1};q;q,q).
\eea

We can obtain another expression for $\rho_0$. We apply the transformation 
\cite[(III.5)]{Ga:Ra}
\bea
{}_2\phi_1(a, b; c; q, z) = \frac{(abz/c;q)_\infty}{(bz/c;q)_\infty}
 {}_3\phi_2(a, c/b, 0; c, cq/bz;q, q)
\eea
to see that
\bea
\sum_{n=k}^\infty \frac{q^n}{(q;q)_n}\binom{n}{k}_q^2=\frac{1}{(q;q)\sf}
\sum_{j=0}^k\frac{q^{k+j}}{(q;q)_j^2}.
\eea

We then find
\bea
\rho_0=\frac{1}{(q;q)\sf}
\sum_{k=0}^\infty q^{2(k+\frac12)^2}\sum_{j=0}^k \frac{q^j}{(q;q)_j^2}.
\eea

A formula more general than (3.6) is
  \bea
  \sum_{n=k}^\infty \frac{\w^n}{(q;q)_n}\binom{n}{k}_q^2=\frac{1}{(\w;q)\sf}
  \sum_{j=0}^k\frac{(\w;q)_j\w^{2k-j}}{(q;q)_j(q;q)_{k-j}^2} \nonumber
  \eea
and is stated in \cite{Al:Ca}. This more general identity
also follows from (3.5) and the simple observation
\bea
\frac{(q^{-k};q)_j}{(q^{1-k}/\w;q)_j} = \frac{(q;q)_k\, (\w;q)_{k-j}}
{(\w;q)_k\, (q;q)_{k-j}} (\w/q)^j.
\nonumber
\eea
We have numerically computed the smallest eigenvalue of the Hankel matrix 
of various dimensions with the Stieltjes--Wigert weight from which
we extrapolate to determine the smallest eigenvalue
 $s=\lim_{N\to\infty}\la_N$ of the infinite Hankel
matrix for different values of $q.$
This is then compared with the numerically
computed lower bound $l=1/\rho_0.$
For $q=\frac12$ we have $s=0.3605\ldots, l=0.3435\ldots$.
 The percentage error $100(s-l)/s$ is
plotted for various values of $q$ and is shown in figure 1.

\begin{figure}[htb]
\begin{center}
\includegraphics[width=7cm,angle=270]{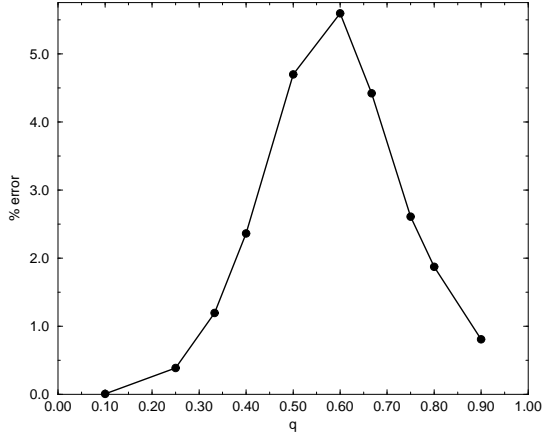}
\end{center}
\caption{Percentage error plotted for various values of $q.$ \label{fig.1}}
\end{figure}

\bigskip
{\bf Example 2}. Al-Salam--Carlitz polynomials.

The Al-Salam-Carlitz polynomials were introduced in \cite{Al:Ca}. We consider
the indeterminate polynomials $V_n^{(a)}(x;q)$, where $0<q<1$ and $q<a<1/q$,
cf. \cite{Be:Va}. For the corresponding orthonormal polynomials $\{p_k\}$
we have by \cite[(4.24)]{Be:Va}
\bea
\sum_{k=0}^{\infty}|p_k(\ei)|^2=\frac{(q\ei,q\emi;q)\sf}{(aq,q,q;q)\sf}
\;_3\phi_2(\ei,\emi,aq; q\ei,q\emi;q,q/a).
\eea
Therefore
\bea
\rho_0=\int_{0}^{2\pi} \sum_{k=0}^{\infty}|p_k(\ei)|^2\dt
=\frac{1}{(aq,q,q;q)\sf}\sum_{n=0}^{\infty}I_n\frac{(aq;q)_n}{(q;q)_n}
\left(\frac{q}{a}\right)^n,
\eea
where
\bea
\mbox{} I_n&=&\int_{0}^{2\pi}\frac{(\ei,\emi;q)\sf}{(1-q^n\ei)(1-q^n\emi)}\dt\\
\mbox{} &=&\int_{|z|=1}\frac{(z,1/z;q)\sf}{(1-q^nz)(1-q^n/z)}\frac{dz}{2\pi iz}.
\nonumber
\eea

Recall the Jacobi  triple product identity \cite{Ga:Ra},
\bea
j(z):=(q,z,1/z;q)\sf=\sum_{k=-\infty}^{\infty} c_k z^k,
\eea
with
\bea
c_k=(-1)^k \left[q^{k(k+1)/2} + q^{k(k-1)/2}\right].
\eea
Note that $c_k=c_{-k}$.

Using the partial fraction decomposition
\bea
\frac{q^n}{1-q^nz}-\frac{q^{-n}}{1-q^{-n}z}
=\frac{1-q^{2n}}{(1-q^nz)(z-q^n)}\nonumber
\eea
we find by the residue theorem and the Jacobi 
triple product identity (3.11) that for $n\geq 1$, $I_n$ is given by
\bea
\mbox{} (1-q^{2n})(q;q)\sf\;I_n&=&
q^n{\rm Res}\left(\frac{j(z)}{1-q^nz},z=0\right)-q^{-n}{\rm Res}
\left(\frac{j(z)}{1-q^{-n}z},z=0\right)\nonumber\\
\mbox{} &=&q^n\sum_{k=0}^{\infty}q^{nk}c_{-k-1}-q^{-n}\sum_{k=0}^{\infty}
q^{-nk}c_{-k-1}\nonumber\\
\mbox{} &=&\sum_{k=1}^{\infty}\left(q^{nk}-q^{-nk}\right)c_k,\nonumber
\eea
while for $n=0$, $I_0$ is
\bea
\mbox{} (q;q)\sf I_0 &=& \int_{|z|=1}\frac{j(z)}{(1-z)(z-1)}\frac{dz}{2\pi i}
=-{\rm Res}\left(\frac{j(z)}{(1-z)^2},z=0\right)\nonumber\\
\mbox{} &=&-\sum_{k=0}^{\infty}
(k+1)c_{-k-1}=\sum_{k=0}^{\infty}(-1)^kq^{k(k+1)/2}.
\nonumber
\eea
The conclusion is
\bea
\mbox{} I_0&=&\frac{1}{(q;q)\sf}\sum_{k=0}^{\infty}(-1)^k
q^{k(k+1)/2},\\
\mbox{} I_n&=&\frac{1}{(1-q^{2n})(q;q)\sf}\sum_{k=1}^{\infty}c_k\left(q^{nk}-
q^{-nk}\right),\quad n\geq 1\nonumber.
\eea
The above formulas can be further simplified. Using the Jacobi triple 
product identity  (3.11) we find for integer values of $n$
\bea
\sum_{k=-\infty}^\infty (-1)^kq^{nk}q^{\binom{k}{2}} =0, \nonumber
\eea
hence
\bea
\sum_{k = 0}^\infty (-1)^kq^{nk}\, q^{\binom{k}{2}} = - 
\sum_{k = 1}^\infty(-1)^kq^{-nk}\, q^{\binom{k+1}{2}}, \quad n = 0, \pm1, \dots. 
\eea
This analysis implies 
\bea
(q;q)_\infty(1-q^{2n})\,  I_n = 2 \sum_{k=1}^\infty (-1)^k\, q^{\binom{k}{2}} 
\left[q^{nk}- q^{-nk}\right]. 
\eea
Thus we have established the representation for $n\geq 1$
\bea
I_n = \frac{2q^{-n}}{(q;q)_\infty} \sum_{k=1}^\infty (-1)^{k-1}\,
q^{\binom{k}{2}}
\; \frac{\sin(nk\tau)}{\sin (n\tau)},       \quad q= e^{-i\tau}.
\eea
It is clear that $I_0$ is the limiting case of $I_n$ as $n \to 0$. The
representation (3.16) indicates that  $I_n$ is a theta function evaluated at
the special point $n\tau$, hence we do not expect to find a closed form
expression for $I_n$.

{\bf Example 3}. Freud-like weight.

In \cite{Be:Va} Berg-Valent found the Nevanlinna matrix in the case of
the indeterminate moment problem corresponding to a birth and death
process with quartic rates. Later Chen and Ismail, cf. \cite{Ch:Is},
considered the corresponding  symmetrized moment problem, found
the Nevanlinna matrix and observed that there are solutions which behave
as the Freud weight $\exp(-\sqrt{|x|})$. In particular they found
the entire functions
\bea
B(z)=-\delta_0(K_0{\sqrt {z/2}}),\quad
D(z)=\frac{4}{\pi}\delta_2(K_0{\sqrt {z/2}}),
\eea
where
\bea
\mbox{} \delta_l(z)&=&\sum_{n=0}^{\infty}\frac{(-1)^n}{(4n+l)!}z^{4n+l},\;\;\;
l=0,1,2,3,\\
\mbox{} K_0&=&\frac{\Ga(1/4)\Ga(5/4)}{{\sqrt \pi}}.
\eea
Note that
\bea
\mbox{} \delta_0(z)&=&\frac{1}{2}\left[\cosh(z{\sqrt i})+
\cos(z{\sqrt i})\right],\\
\mbox{} \delta_2(z)&=&\frac{1}{2i}\left[\cosh(z{\sqrt i})-\cos(z{\sqrt
i})\right].
\eea
If $\omega:={\rm exp}(i\pi/4)=(1+i)/{\sqrt 2},$ then a simple
calculation shows that
\bea
\mbox{} &{}& B(x)D(y)-D(x)B(y) \\
\mbox{} &{}& \quad = \frac{-2i}{\pi}\left[\cos (\w^3K_0\sqrt{x/2})
\cos (\w K_0\sqrt{y/2})- \cos (\w^3 K_0\sqrt{y/2})\cos (\w K_0\sqrt{x/2})
\right]. \nonumber
\eea
If $x=\ei,$ and $y=\emi,$ then we linearise the products of cosines and find
that the right-hand side of (3.22) is
\bea
&{}& \frac{-i}{\pi}\left\{\cos [K_0(\w^3e^{i\t/2} + \w e^{-i\t/2})/\sqrt{2}]
+  \cos [K_0(\w^3e^{i\t/2}- \w e^{-i\t/2})/\sqrt{2}] \right. \nonumber\\
&{}& \quad \left. - \cos [K_0(\w^3e^{-i\t/2} + \w e^{i\t/2})/\sqrt{2}]
- \cos [K_0(\w^3e^{-i\t/2} - \w e^{i\t/2})/\sqrt{2}]\right\}
 \nonumber
\eea
We now combine the first and third terms, then combine the second and fourth
terms and apply the addition theorem for trigonometric functions. We then
see that the above is
\bea
&{}& \frac{2i}{\pi}\left\{\sinh [K_0\cos (\t/2)] \sinh [K_0\sin (\t/2)]
+ \sin  [K_0\cos (\t/2)] \sin  [K_0\sin (\t/2)]\right\}. \nonumber
\eea
Thus we have proved that
\bea
\mbox{} &{}& \frac{B\left(\ei\right)D\left(\emi\right)-
B\left(\emi\right)D\left(\ei\right)}
{\ei-\emi} \\
\mbox{} &{}&
\quad = \frac{1}{\pi \sin \theta}
\,\left\{\sinh [K_0\cos(\theta/2)] \sinh [K_0\sin(\theta/2)]
+  \sin  [K_0\cos(\theta/2)] \sin  [K_0\sin(\theta/2)]\right\}. \nonumber
\eea
Thus in the case under consideration, after some straightforward calculations
and the evaluation of a beta integral, we obtain
\bea
\mbox{} \rho_0=\int_{0}^{2\pi}\sum_{n=0}^{\infty}\left|p_n(\ei)\right|^2\,
\frac{d\theta}{2\pi}&=& \frac{K_0^2}{\pi}
\sum_{m, n \ge 0, m+n \, {\rm even}}
\frac{(K_0/2)^{2m+2n}}{(2m+1)(2n+1)\, m!\, n!\, (m+n)!}.
\eea

{\bf Example 4}. $q^{-1}$-Hermite polynomials.

 Ismail and Masson \cite{Is:Ma}
proved that for this moment problem the functions $B$ and $D$ are given by
\bea
B(\sinh \xi) = -\frac{(qe^{2\xi}, qe^{-2\xi};q^2)_\infty}{(q, q;q^2)_\infty},
\quad D(\sinh \xi) = \frac{\sinh \xi}{(q;q)_\infty}\,
(q^2e^{2\xi}, q^2e^{-2\xi};q^2)_\infty,
\eea
\cite[(5.32)]{Is:Ma}, \cite[(5.36)]{Is:Ma}; respectively.
Ismail and Masson also showed that \cite[(6.25)]{Is:Ma}
\bea
&{}& B(\sinh \xi)D(\sinh \eta)- B(\sinh \eta) D(\sinh \xi) \\
&{}& \quad = \frac{-e^{ \eta}}{2(q;q)_\infty} \, \prod_{n=0}^\infty
\left[1-2e^{-\eta} q^n \sinh \xi -e^{-2\eta}  q^{2n}\right]
\left[1+2e^{\eta} q^{n+1} \sinh \xi -e^{2\eta}  q^{2n+2}\right]. \nonumber
\eea
We rewrite the infinite product as
\bea
\prod_{n=0}^\infty a_nb_n=a_0\prod_{n=1}^\infty a_nb_{n-1},\nonumber
\eea
and with
$\sinh\xi = e^{i\t}$ and $\sinh \eta = e^{-i\t}$ we get the following
representation
\bea
&{}& \frac{B(e^{i\t})D(e^{-i\t})- B(e^{-i\t})D(e^{i\t})}{e^{i\t}-e^{-i\t}} \\
&{}& \quad = \frac{1}{(q;q)_\infty} \prod_{n=1}^\infty
\left[1+4 q^n -2q^{2n} + 4q^{3n} + q^{4n} - 8q^{2n} \cos (2\t)\right]
\nonumber\\
&{}& \quad = \frac{1}{(q;q)_\infty} \prod_{n=1}^\infty
\left[(1+q^n)^4-16q^{2n}\cos^2\t\right].\nonumber
\eea

Writing the infinite product as a power series in $\cos^2\t$ and using
\bea
\int_{-\pi}^{\pi} \cos^{2k}\t\;\frac{d\t}{2\pi}=2^{-2k}\binom{2k}{k},\nonumber
\eea
we evaluate the integral of (3.27) with respect to $d\t/2\pi$ as
\bea
\rho_0=\frac{(-q;q)_\infty^4}{(q;q)_\infty}\sum_{k=0}^\infty \binom{2k}{k}
\sum_{1\leq n_1<\ldots<n_k}\frac{(-2)^{2k}q^{2(n_1+\cdots +n_k)}}
{\left[(1+q^{n_1})\cdots (1+q^{n_k})\right]^4}.
\eea

The formula (3.27) can be transformed further by putting
$\cos^2\psi=-\cos\t$ and $p^2=q$, because then
$$
\prod_{n=1}^\infty\left[(1+q^n)^2+4q^n\cos\t\right]=
\prod_{n=1}^\infty\left[1+p^{4n}-2p^{2n}\cos(2\psi)\right]
$$
can be expressed by means of the theta function $\vartheta_1(p;\psi)$.
We find
\bea
\prod_{n=1}^\infty\left[(1+q^n)^2+4q^n\cos\t\right]=
\frac{1}{(q;q)_\infty}\sum_{n=0}^\infty (-1)^nq^{\binom{n+1}{2}}U_{2n}
(\cos\psi),
\eea
where
$$
U_{2n}(\cos\psi)=\frac{\sin(2n+1)\psi}{\sin\psi}
$$
is the Chebyshev polynomial of the second kind given by
\bea
U_{2n}(x)=\sum_{k=0}^n\binom{2n+1}{2k+1} (-1)^k x^{2(n-k)}(1-x^2)^k.
\eea

Similarly putting $\cos^2\varphi=\cos\t$ we find
\bea
\prod_{n=1}^\infty\left[(1+q^n)^2-4q^n\cos\t\right]=
\frac{1}{(q;q)_\infty}\sum_{n=0}^\infty (-1)^nq^{\binom{n+1}{2}}U_{2n}
(\cos\varphi).
\eea
If we let $U_n^*$ be the polynomial of degree $n$ such that
 $U_{2n}(x)=U_n^*(x^2)$, we get

\bea
&{}&
\frac{B(e^{i\t})D(e^{-i\t})- B(e^{-i\t})D(e^{i\t})}{e^{i\t}-e^{-i\t}}\\
&{}&\quad=\frac{1}{(q;q)_\infty^2}\sum_{n,m=0}^\infty (-1)^mq^{\binom
{n+1}{2}+\binom{m+1}{2}}U_n^*(-\cos\t)U_m^*(\cos\t).\nonumber
\eea

For non-negative integers $k,l,r$ we have
\bea
&{}& C(k,l,r):=\frac{1}{2\pi}\int_0^{2\pi}(1+\cos\t)^k(1-\cos\t)^l
\cos^r\t\;d\t\\
&{}& \quad =\frac{2^{k+l}}{\pi}(-1)^rB(k+\frac12,l+\frac12)
{}_2F_1(k+\frac12,-r;k+l+1;2),\nonumber
\eea
which gives
\bea
&{}& \frac{1}{2\pi}\int_0^{2\pi} U_n^*(-\cos\t)U_m^*(\cos\t)\;d\t\\
&{}&=\sum_{k=0}^n\sum_{l=0}^m \binom{2n+1}{2k+1}\binom{2m+1}{2l+1}(-1)^{n+l}
C(k,l,n+m-k-l).\nonumber
\eea

Putting these formulas together we get a 5-fold sum for $\rho_0$.

\medskip
{\bf Acknowledgement} The authors would like to thank Mr. N. D. Lawrence
for supplying the numerical data and the graph.

\noindent Department of Mathematics, University of Copenhagen,
Universitetsparken 5, DK-2100 Copenhagen, Denmark

\noindent Department of Mathematics, Imperial College,
180 Queen's Gate, SW7 2BZ, London, England.

\noindent Department of Mathematics, University of South Florida, Tampa,
 Florida 33620-5700.

\noindent
e-mail \qquad berg@math.ku.dk, \qquad y.chen@ic.ac.uk \qquad and
\qquad ismail@math.usf.edu

\end{document}